\documentclass[a4paper,12pt]{article}
\usepackage{amssymb,latexsym,amsmath,amsthm,amscd,graphpap}
\usepackage{amsmath}
\usepackage{amsfonts}
\usepackage{amssymb}
\usepackage[dvips]{epsfig}
\usepackage{epsfig}

\newtheorem{theorem}{Theorem}[section]

\newtheorem{prop}[theorem]{Proposition}

\newenvironment{demo}{ \noindent \emph{\textbf{Proof: }}}{\hfill$\square$\\
\vspace{0.4cm}}

\newcommand{\Rm}{\mathbb{R}}

\newcommand{\Lm}{\mathbb{L}}

\newcommand{\Nm}{\mathbb{N}}

\newcommand{\Tm}{\mathbb{T}}

\newcommand{\Cm}{\mathcal{C}}
\newcommand{\Nc}{\mathcal{N}}
\newcommand{\Mc}{\mathcal{M}}

\newcommand{\Am}{\mathcal{A}}

\newcommand{\Um}{\mathcal{U}}

\newcommand{\Drond}[2]{\frac{\partial #1}{\partial #2}}
\newcommand{\grad}{{\nabla}}

\newcommand{\no}{n$^{\text{o}}$}

\numberwithin{equation}{section}

\newcounter{constantes}

\newdimen\texpscorrection
\newdimen\figcenter
\def\figurewithtex #1 #2 #3 #4 #5\cr{\null
   {\goodbreak\figcenter=\hsize\relax
   \advance\figcenter by -#4truecm
   \divide\figcenter by 2
   \begin{figure}[hbt]
   \vskip #3truecm\noindent\hskip\figcenter
   \includegraphics{#1}{\hskip\texpscorrection\input #2 }
   \vskip 0.8truecm{\baselineskip=0.8\baselineskip
   \noindent \vbox{\noindent {\footnotesize #5}}\par}
   \end{figure}}}
\def\point#1 #2 #3 {\rlap{\kern #1 truecm
\raise #2 truecm \hbox{#3}}}

\newcommand{\pc}{ \usefont{T1}{cmtl}{m}{n} \selectfont}

\begin{document}

\title{\bf A striking correspondence between the dynamics generated by the
vector fields and by the scalar parabolic equations}

\author{Romain JOLY\\
{\small Univ. Grenoble Alpes, CNRS}\\
{\small Institut Fourier}\\
{\small F-38000 Grenoble, France}\\ 
{\pc \small Romain.Joly@univ-grenoble-alpes.fr} \\[5mm]
Genevi\`eve RAUGEL\\
{\small CNRS, Univ Paris-Sud}\\
{\small Laboratoire de Math\'{e}matiques d'Orsay }\\
{\small F-91405 Orsay cedex, France}\\
{\pc \small Genevieve.Raugel@math.u-psud.fr}
}

\date{}

\maketitle

\vspace{3cm}

\hrule

{\quote
{
{\em \`A la m\'{e}moire de  Michelle Schatzman, notre ch\`{e}re et regrett\'{e}e
coll\`{e}gue et amie.}\\[3mm]
{\small
The results of this paper have been presented 
at the Conference
in honour of the sixthieth birthday of Michelle
Schatzman.
The second 
author would like to express her gratitude towards Michelle, who guided her 
first steps in the research in mathematics and was a true friend.}}

}

\vspace{1cm}

\hrule

\pagebreak

\begin{abstract}
The purpose of this paper is to enhance a correspondence between the
dynamics
of the differential equations $\dot y(t)=g(y(t))$ on $\Rm^d$ and those of
the parabolic equations $\dot u=\Delta u +f(x,u,\grad u)$ on a bounded domain
$\Omega$.
We give details on the similarities of these dynamics in the cases
$d=1$, $d=2$ and $d\geq 3$ and in the corresponding cases
$\Omega=(0,1)$, $\Omega=\Tm^1$ and dim($\Omega$)$\geq 2$ respectively.
In addition to the beauty of such a correspondence, this could serve as a
guideline for future research on the dynamics of parabolic
equations.
\end{abstract}

\vspace{2cm}


\section{Introduction}
In this paper, we want to point out the similarities between the dynamics
of vector fields in $\Rm^d$ and those of reaction-diffusion equations on
bounded domains. More precisely, we consider the following classes of 
equations.\\[2mm]

{\noindent \bf Class of vector fields}\\
Let $d\geq 1$ and $r\geq 1$ and let $g\in\Cm^r(\Rm^d,\Rm^d)$ be a given vector
field. We consider the ordinary differential equation
\begin{equation}\label{ODE}
\left\{\begin{array}{l}\dot y(t)=g(y(t))~~t>0\\ y(0)=y_0\in\Rm^d
\end{array}\right.
\end{equation}
where $\dot y(t)$ denotes the time-derivative of $y(t)$.\\
The equation \eqref{ODE} defines a local dynamical system $T_g(t)$ on 
$\Rm^d$ by
setting $T_g(t)y_0=y(t)$. We assume that there exists $M>0$ large enough such
that 
$$\forall y\in\Rm^d,~\|y\|\geq M~\Rightarrow~\langle y|g(y)\rangle <0~.$$ 
This condition ensures that $T_g(t)$ is a global dynamical system. Moreover, the
ball $B(0,M)$ attracts the bounded sets of $\Rm^d$. Therefore,
$T_g(t)$ admits a compact global attractor\footnotemark[1] \footnotetext[1]{ To
make the reading of this article easier for the reader, who is not
familiar with dynamical systems theory or with the study of PDEs,
we add a short glossary at the end of the paper.}
$\Am_g$. The attractor $\Am_g$ contains the most interesting
trajectories such as
periodic, homoclinic and heteroclinic orbits\footnotemark[1] and any $\alpha-$
or $\omega-$limit
set\footnotemark[1]. Therefore, if one neglects the transient  
dynamics, the dynamics on $\Am_g$ is a good representation of the whole dynamics
of $T_g(t)$.\\[2mm]

{\noindent \bf Class of scalar parabolic equations}\\
Let $d'\geq 1$ and let $\Omega$ be either a regular  bounded
domain
of $\Rm^{d'}$, or the torus $\Tm^{d'}$. We choose $p> d'$ and $\alpha\in
((p+d')/2p,1)$. We denote
$X^\alpha \equiv D((-\Delta_N)^{\alpha})$
the fractional power space\footnotemark[1]
associated with the Laplacian operator $\Delta_N$ on $\Lm^p(\Omega)$
with homogeneous Neumann boundary conditions. It is well-known\footnotemark[1]
that $X^\alpha$ is continuously embedded in the Sobolev space
$W^{2\alpha,p}(\Omega)$  and thus it is compactly embedded in
$\Cm^1(\overline\Omega)$. Let $r\geq 1$ and
$f\in\Cm^r(\overline\Omega \times \Rm \times \Rm^{d'}, \Rm)$.
We consider the parabolic equation
\begin{equation}\label{PDE}
\left\{\begin{array}{ll}\dot u(x,t)=\Delta u(x,t) +f(x,u(x,t),\grad
u(x,t))~&~(x,t)\in\Omega\times (0,+\infty)\\ \Drond u\nu (x,t)=0~&~(x,t)
\in\partial\Omega\times(0,+\infty)\\ u(x,0)=u_0(x)\in X^\alpha &
\end{array}\right.
\end{equation}
where $\dot u(t)$ is the time-derivative of $u(t)$.\\
Eq. \eqref{PDE} defines a local dynamical system $S_f(t)$ on
$X^\alpha$ (see \cite{Henry-book}) by setting
$S_f(t)u_0=u(t)$. We assume moreover that there exist $c\in\Cm^0(\Rm_+,\Rm_+)$,
$ \varepsilon>0$ and $\kappa>0$ such that $f$ satisfies
\begin{align*}
\forall R>0,~\forall 
\xi\in\Rm^{d'},~&\sup_{(x,z)\in\overline\Omega\times[-R,R]}
|f(x,z,\xi)|\leq c(R)(1+|\xi|^{2-\varepsilon}) \\
\text{and } \forall z\in\Rm ,~\forall
x\in\overline\Omega,~~~&~~|z|\geq\kappa~\Rightarrow~
zf(x,z,0) < 0~.
\end{align*}
Then, Eq. \eqref{PDE} defines a global dynamical system in
$X^\alpha$ which admits a compact global attractor $\Am$ (see
\cite{Polacik-handbook}).\\[2mm]

\begin{table}\label{table}
{\noindent \small \begin{tabular}{|c|l|c|}
\hline
ODE& &PDE\\
\hline \hline
\begin{tabular}{c}$d=1$\\$~$\\{\footnotesize (Or more 
generally}\\{\footnotesize
tridiagonal}\\{ \footnotesize cooperative}\\{ \footnotesize
system of ODEs)} \end{tabular} &\begin{tabular}{l}
$\bullet$ Gradient  dynamics\\ $\bullet$
  Convergence to
an equilibrium point\\ $\bullet$ Automatic transversality of stable\\ and
unstable manifolds\\
$\bullet$ Genericity of Morse-Smale property\\
$\bullet$ Knowledge of the equilibrium points\\ implies knowledge of the whole
dynamics\\
$\bullet$ Dimension of the attractor equal to the\\ largest dimension of the
unstable manifolds\\
$\bullet$ Realisation of the ODE in the PDE \end{tabular} &$\Omega=(0,1)$\\
\hline \hline
\begin{tabular}{c}$d=2$\\General case\\$~$\\{\footnotesize(Or more
generally}\\{\footnotesize cyclic tridiagonal}\\{ \footnotesize
cooperative }\\{\footnotesize system of ODEs)}
\end{tabular}
&
\begin{tabular}{l}
$\bullet$ Poincar\'e-Bendixson property\\
$\bullet$ Automatic transversality of stable\\ and
unstable manifolds of two orbits\\ if one of them is a hyperbolic
periodic orbit\\ or if both are equilibrium points\\ with different Morse
indices.\\
$\bullet$ Non-existence of homoclinic orbits\\ for periodic orbits\\
$\bullet$ Genericity of Morse-Smale property\\
$\bullet$ Realisation of the ODE in the PDE
\end{tabular} &
\begin{tabular}{c}$\Omega=\Tm^1$\\General case\end{tabular}\\
\hline
\begin{tabular}{c}$d=2$\\$g$ radially\\ symmetric \end{tabular} &
\begin{tabular}{l}
$\bullet$ Automatic transversality of stable\\ and
unstable manifolds of\\ equilibrium points and periodic orbits\\
$\bullet$ No homoclinic orbit\\
$\bullet$ Knowledge of the equilibrium points\\ and of the periodic
orbits\\ implies knowledge of the whole
dynamics\\
$\bullet$ Genericity of the Morse-Smale property.\\
$\bullet$ Dimension of the attractor equal to the\\ largest dimension of the
unstable manifolds\\
$\bullet$ Realisation of the ODE in the PDE
\end{tabular} &
\begin{tabular}{c}$\Omega=\Tm^1$\\ $f(x,u,\grad u)\equiv f(u,\grad u)$
\end{tabular}\\
\hline  \hline
$d\geq 3$&
\begin{tabular}{l}
$\bullet$ Existence of persistent chaotic dynamics\\
$\bullet$ Genericity of Kupka-Smale property (ODE)\\
$\bullet$ Generic transversality of homoclinic\\
and heteroclinic orbits (PDE)\\
$\bullet$ Realisation of the ODE in the PDE
\end{tabular} &
dim($\Omega$)$\geq 2$\\
\hline \hline
\begin{tabular}{c}Any $d$\\$g\equiv \grad G$\end{tabular} &
\begin{tabular}{l}
$\bullet$ Gradient  dynamics\\
$\bullet$ Genericity of the Morse-Smale property\\
$\bullet$ Realisation of a generic ODE in the PDE
\end{tabular} &
\begin{tabular}{c}Any $\Omega$\\$f(x,u,\grad u)\equiv f(x,u)$ \end{tabular}\\
\hline
\end{tabular}

}\caption{ the correspondence between the dynamics of
vector fields and the ones of parabolic equations.}\end{table}

The reader, which is not familiar with partial differential equations, may
neglect all the technicalities about $X^\alpha$, the Sobolev spaces and the
parabolic equations in a first reading.
The most important point is that $S_f(t)$ is a dynamical system
defined on an {\it infinite-dimensional function space}. Compared with the
finite-dimensional case, new difficulties arise. For example, the 
existence of a
compact global attractor requires compactness properties, coming here from the
smoothing effect of \eqref{PDE}. We also mention that, even if the backward
uniqueness property holds, backward trajectories do not exist in general for
\eqref{PDE}. 
The reader interested in the dynamics of
\eqref{PDE} may consult \cite{Fiedler-Scheel}, \cite{Henry-book},
\cite{HaleCanada}, \cite{Polacik-handbook} or \cite{GR-handbook}.

The purpose of this paper is to emphasize the different relationships
between the dynamics of \eqref{ODE} and \eqref{PDE}. The correspondence is
surprisingly perfect. It can be summarized by Table \ref{table}.
This correspondence has already been noticed for some of the properties of the
table. We complete here the correspondence for all the known
properties of the dynamics of the parabolic equation. Table
\ref{table} will be discussed in more details in Section
\ref{section-comments} and, for cooperative
systems, in Section \ref{section-coop}. Some of the properties
presented in the table concerning finite-dimensional dynamical systems are
trivial, other ones are now well-known. However, the corresponding results
for the parabolic equation are more involved and some of them are very recent.
These properties are mainly based on Sturm-Liouville arguments and unique
continuation properties for the
parabolic equations as explained in Section \ref{section-ucp}. The study of the
dynamics generated by vector fields in dimension $d\geq 3$ is still a subject
of research. Taking into account the correspondence presented in Table
\ref{table}
should give a guideline for research on the dynamics of the
parabolic equations. Some examples of open questions are given in Section
\ref{section-open}.

We underline
 that we only consider the dynamics on the compact global attractors.
Hence, we deal with dynamical systems on compact sets. It is important to be
aware of the fact that, even if the dimension of the compact global
attractor $\Am$ of the parabolic equation \eqref{PDE} is finite, it can be made
as large as wanted by choosing a suitable function $f$. This is true even
if $\Omega$ is one-dimensional. Therefore, all the possible properties of the
dynamics of \eqref{PDE} do not come from the low dimension of $\Am$ but from
properties, which are very particular to the flow of the parabolic equations.

Finally, we remark that most of the results described here also hold in more
general frames than \eqref{ODE} and \eqref{PDE}. For example, $\Rm^d$ could
be replaced by a compact orientable manifold without boundary. We
could also choose for \eqref{PDE} more general boundary conditions than Neumann
ones, or less restrictive growing conditions for $f$. The domain $\Omega$
may be replaced by a bounded smooth manifold. Finally, notice that the case $\Omega=\Tm^{d'}$
can be seen as $\Omega=(0,1)^{d'}$ with periodic boundary conditions.


\section{Details and comments about the
correspondence table}\label{section-comments}

We expect the reader to be familiar with the basic notions of the
theory of dynamical systems and flows. Some definitions are briefly recalled in
the glossary at the end of this paper. For more precisions, we refer 
for example
to \cite{Katok-Hasselblatt}, \cite{Newhouse2}, \cite{Palis-de-Melo},
\cite{Robinson2} or \cite{Ruelle} for finite-dimensional dynamics and to
\cite{HMO},
\cite{Henry-book}, \cite{Robinson} or \cite{HJR} for the infinite-dimensional
ones.

We first would like to give short comments and motivations
concerning the properties appearing in Table \ref{table}. Notice that we do not
deal in this section with the cooperative systems of ODEs. The properties of
these systems are discussed in Section \ref{section-coop}.

A {\it generic property of the dynamics} is a property satisfied by a
countable intersection of open dense subsets of the considered class of
dynamical systems. Generic dynamics represent the typical behaviour of a class
of dynamical systems.
For finite-dimensional flows, we mainly consider classes of the form
$(T_g(t))_{g\in\Cm^1(\Rm^d,\Rm^d)}$. The parameter is
the vector field $g$, 
which belongs to the space $\Cm^1(\Rm^d,\Rm^d)$ endowed with either the 
classical $\Cm^1$
or the $\Cm^1$ Whitney topology. 
Notice that the question wether or not a property is
generic for $g\in\Cm^r(\Rm^d,\Rm^d)$ for some $r\geq 2$ may be much more
difficult than $\Cm^1$ genericity. We will not discuss this problem here. In
some cases, we restrict the class of vector fields to subspaces of
$\Cm^1(\Rm^d,\Rm^d)$ such as radially symmetric, gradient vector fields or
cooperative systems. In a
similar way, for infinite-dimensional dynamics, we consider families of the
type $(S_f(t))_{f\in\Cm^1(\overline\Omega\times\Rm\times\Rm^{d'},\Rm)}$, where
$\Cm^1(\overline\Omega\times\Rm\times\Rm^{d'},\Rm)$ is endowed with either
the classical $\Cm^1$ or the $\Cm^1$ Whitney topology.
 For some results, we restrict the class
of nonlinearities $f$ to homogeneous ones or to ones, which are independent of the last
variable $\xi$.

{\it Poincar\'e-Bendixson property} and {\it the convergence to an
equilibrium or a periodic orbit} are properties related to the
following question: how
simple are the $\alpha-$ and $\omega-$limit sets of the trajectories~? For
vector
fields, the restriction of the complexity of the limit sets may come from the
restriction of freedom due to the low dimension of the flow. As said above,
there is no restriction on the dimension of the global attractor for the
parabolic equations. The possible restrictions of the complexity
of the limit sets come from particular properties of the
parabolic equations, see Section \ref{section-ucp}.

{\it Hyperbolicity of equilibria and periodic orbits, tranversality of stable
and unstable
manifolds, Kupka-Smale and Morse-Smale properties} are properties related to
the question of  stability of the local and global dynamics respectively.
Morse-Smale property
is the strongest one. It implies the structural stability of the global
dynamics: if the
dynamical system $T_g(t)$
satisfies the Morse-Smale property, then for $\tilde g$ close enough to $g$, the
dynamics of $T_{\tilde g}(t)$, restricted to its attractor $\Am_{\tilde g}$, are
qualitatively the same as the ones of $T_g(t)$ on $\Am_g$, see \cite{Palis},
\cite{Palis-Smale} and
\cite{Palis-de-Melo}. The same structural stability result holds for parabolic
equations satisfying the Morse-Smale property, see \cite{HMO}, \cite{HJR} and
\cite{Oliva}. It is natural to
wonder if almost all the dynamics satisfy these properties, that is if these
properties are generic.

The fact that {\it the knowledge of the equilibria and the periodic
 orbits implies the knowledge of the whole dynamics} may be studied at different
 levels. Two equilibria or periodic orbits being given, can we know if they are
 connected or not by a heteroclinic orbit ? Are two dynamics with the same
 equilibria and periodic orbits equivalent ? Is there a simple algorithm to
 determine the global dynamics from the position of the equilibria and the
 periodic orbits ? These questions are among the rare dynamical
 questions coming from the study of partial differential equations and not from
 the study of vector
fields. Indeed, for finite-dimensional dynamical systems, the answers,
either positive or negative, are too simple. In contrast, such kinds of
results are probably among  the most amazing ones for the dynamics of 
the parabolic equations.

{\it The persistent
chaotic dynamics} and the
fact that {\it the dimension of the attractor is equal to the largest
dimension of the unstable manifolds}, are
related to the following question: how complicated may be the dynamics~?
In general, the dimension of the attractor of a dynamical system may
be larger than the largest dimension of the unstable manifolds. The classes of
systems, where these dimensions automatically coincide,  are
strongly constrained, which in some sense implies a simple behaviour.
On the contrary, chaotic dynamics have very complicated behaviour. Chaotic
dynamics may occur through several phenomena, and the notion of chaotic
behaviour depends on the authors. In this paper, ``persistent 
chaotic dynamics'' refers to the presence of a tranversal homoclinic orbit generating
a Smale horseshoe (see \cite{Smale65}). The persistent chaotic dynamics
provide complicated dynamics, which cannot
be removed by small perturbations of the system. Such an open set of
chaotic dynamics is a counter-example to the genericity of the Morse-Smale systems.

The question of {\it the realization of vector fields in the
parabolic equations} is as follows:
a vector field $g\in\Cm^r(\Rm^d,\Rm^d)$ being given, can we find a function $f$ and an invariant
manifold $M\subset L^{p}(\Omega)$ such that the dynamics of the parabolic
equation \eqref{PDE} restricted to $M$ is equivalent to the dynamics generated
by the vector field $g$~? A positive answer to this question implies that
the dynamics of the considered class of parabolic equations is at least as
complicated as the dynamics of the considered class of vector fields. Such a
realization result is very interesting since, on the opposite, the other
properties stated in Table \ref{table} roughly say that
the dynamics of the
parabolic equation \eqref{PDE} cannot be much more complicated than the ones of
the corresponding class of finite-dimensional flows.
One has to keep in mind that the
manifold $M$, on which the finite-dimensional dynamics are realized, is not
necessarily stable with respect to the dynamics of the parabolic
equation. Typically, $M$ cannot be stable if the finite-dimensional
system contains a stable periodic orbit, since all periodic orbits of
\eqref{PDE} are unstable (see for example \cite{Hirsch3}).

\vspace{2mm}

Now, we give short comments and references for the correspondences stated in
Table \ref{table}.

{\noindent \bf $\bullet$ $d=1$ and $\Omega=(0,1)$}

The dynamics generated by a one-dimensional vector field is very
simple. Its attractor consists in equilibrium points and heteroclinic orbits
connecting two of them. The existence of these heteroclinic orbits is easily
deduced from the positions of the equilibrium points. Moreover, these
heteroclinic connections are trivially transversal. Finally,
\eqref{ODE} is clearly a gradient system with associated  Lyapounov functional
$F(y) =-\int_0^y g(s) ds$.
 As a consequence, the Morse-Smale property is equivalent to the hyperbolicity
of all the
equilibrium points, which holds for a generic one-dimensional vector field.

The dynamics $S_f(t)$ generated by \eqref{PDE} for $\Omega=(0,1)$ is richer
since its attractor may have a very large dimension. However, these dynamics
satisfy similar properties. These similarities are mainly due to the
constraints
coming from the non-increase of the number of zeros of solutions of the linear
parabolic equation (see Theorem \ref{th-Angenent}). Zelenyak has proved in
\cite{Zelenyak} that $S_f(t)$ admits an explicit Lyapounov function and thus
that it is gradient. He also showed that the $\omega-$limit sets of the
trajectories consist in single equilibrium points. In Proposition
\ref{prop-cv-single}, we give a short proof of this result, due to Matano. The fact that
the stable and unstable manifolds of equilibrium points always intersect
transversally comes from Theorem \ref{th-Angenent} and the standard
Sturm-Liouville theory. This property has been first proved by Henry in
\cite{Henry-T} and later by Angenent \cite{Angenent-T} in the weaker case of
hyperbolic equilibria. As a consequence of the previous results, the Morse-Smale
property is equivalent to the hyperbolicity of the equilibrium points and is
satisfied by the parabolic equation on $(0,1)$ generically with respect to $f$.
The most surprising result concerning \eqref{PDE} on $\Omega=(0,1)$ is the
following one. Assuming that every equilibrium point is hyperbolic and that the
equilibrium points $e_1$,...,$e_p$ are known, one can say if two given
equilibria $e_i$ and $e_j$ are connected or not by a heteroclinic orbit. This
property has been proved by Brunovsk\'y and Fiedler in \cite{Brunovsky-Fiedler}
for $f=f(u)$ and by Fiedler and Rocha in \cite{Fiedler-Rocha} in the general
case. The description of the heteroclinic connections is
obtained
from the Sturm permutation which is a permutation generated by the respective
positions of the values $e_i(0)$ and $e_i(1)$ of the equilibrium points at the
endpoints of $\Omega=(0,1)$. The importance of Sturm permutation has been
first underlined by Fusco and Rocha in \cite{Fusco-Rocha}. We also refer to the
work of Wolfrum \cite{Wolfrum} , which presents a very nice formalism for this
property.
Fiedler and Rocha  showed in \cite{Fiedler-Rocha-Equivalence}
that the Sturm permutation characterizes the global dynamics of
\eqref{PDE} on $(0,1)$. 
They  proved in \cite{Fiedler-Rocha-realization}
that it is possible to 
give the exact list of  all
the permutations which are Sturm
permutations for some nonlinearity $f$ and thus to give the list of all the possible
dynamics of
the parabolic equation on $(0,1)$. The fact that the dimension of the attractor
is equal to the largest dimension of the unstable manifolds has been shown by
Rocha in \cite{Ro91}. The previous works of Jolly \cite{Jo} and
Brunovsk\'y \cite{Br90} deal with the particular case $f\equiv f(u)$, but show a
stronger result: the attractor can be embedded in a $\Cm^1$ invariant
graph of dimension equal to the largest dimension of the unstable
manifolds. Finally, let us mention that it is easy to
realize any one-dimensional flow in an invariant manifold
of the semi-flow generated by
the one-dimensional parabolic equation. For example, in the simplest case of
Neumann boundary conditions as in \eqref{PDE}, one can realize
the flow of any vector field
$g$ as the restriction of the dynamics of the equation $\dot u=\Delta u +g(u)$ to the subspace of
spatially constant functions.

\vspace{2mm}

{\noindent \bf $\bullet$ $d=2$ and $\Omega=\Tm^1$, general case}\\
Even if they are richer than in the one-dimensional case, the flows generated
by vector fields on $\Rm^2$ are constrained by the
Poincar\'e-Bendixson property (see the original works of Poincar\'e \cite{Poincare} and Bendixson
 \cite{Bendixson} or any textbook on ordinary differential equations).
In particular, this constraint precludes the existence of non-trivial
non-wandering points in Kupka-Smale dynamics. Due to the low dimension of the
dynamics, the stable and unstable manifolds of hyperbolic
equilibria or periodic orbits always  intersect transversally if either one of
the manifold corresponds to a periodic orbit or if the invariant
manifolds correspond to two equilibrium points with different Morse indices.
Moreover, there is no homoclinic trajectory for periodic orbits. Using these
particular properties, Peixoto proved in \cite{Peixoto1} that the Morse-Smale
property holds for a generic two-dimensional vector field.

The first correspondence between two-dimensional flows and the dynamics of the
pa\-ra\-bo\-lic equation \eqref{PDE} on the circle $\Omega=\Tm^1$ has been
obtained by Fiedler and Mallet-Paret in \cite{FMP}. They proved that
the Poincar\'e-Bendixson
property holds for \eqref{PDE} on $\Tm^1$, by using the properties of the zero
number (see Theorem
\ref{th-Angenent}). The realization of any two-dimensional flow in a
two-dimensional invariant manifold of the parabolic equation on the circle has
been proved by Sandstede and Fiedler in \cite{Sandstede-Fiedler}. 
Very recently,
Czaja and Rocha have shown in \cite{Czaja-Rocha} that the stable and unstable
manifolds of two hyperbolic periodic orbits always intersect transversally and
that there
is no homoclinic connection for a periodic orbit. The other automatic
transversality results and the proof of the genericity of the Morse-Smale property
have been completed by the authors in \cite{JR} and \cite{JR2}.

\vspace{2mm}

{\noindent \bf $\bullet$ $d=2$ and $\Omega=\Tm^1$, radial symmetry and
$\Tm^1-$equivariance}\\
When the vector field $g$ satisfies a radial symmetry, the dynamics of the
two-dimensional flow generated by \eqref{ODE} becomes roughly one-dimensional.
The closed orbits consist in $0$, circles of equilibrium points and periodic
orbits being circles described with a constant rotating speed. The
dynamics are so constrained that the closed orbits being given, it is possible
to
describe all the heteroclinic connections. Notice that no homoclinic connection
is possible. We also underline that the Morse-Smale property is generic 
in the class
of radially symmetric vector fields.

If the two-dimensional radial vector fields are too simple to attract much
attention, the corresponding case for the parabolic equation \eqref{PDE} on
$\Omega=\Tm^1$ with homogeneous nonlinearity $f(x,u,\partial_x
u)\equiv f(u,\partial_x u)$ has been extensively studied. Since Theorem
\ref{th-Angenent} holds for \eqref{PDE} with any one-dimensional domain
$\Omega$, it is natural to expect results for
\eqref{PDE} on $\Omega=\Tm^1$ similar to the ones obtained for
\eqref{PDE} on $\Omega=(0,1)$. In particular, one may wonder if it is
possible to describe the global dynamics of \eqref{PDE} knowing
the equilibria and the periodic orbits only. However, this property is
still open for general non-linearities $f(x,u,\partial_x u)$ in the case
$\Omega=\Tm^1$.
Moreover, if one believes in the correspondence stated
in this paper, one can claim that it is in fact false for a general 
nonlinearity
$f(x,u,\partial_x u)$. Therefore, it was natural to first study the simpler
case of homogeneous nonlinearities $f\equiv f(u,\partial_x u)$. Indeed,
the dynamics in this case are
much simpler, in particular the closed orbits are either homogeneous 
equilibrium
points $e(x)\equiv e\in\Rm$, or circles of non-homogeneous equilibrium points,
or periodic orbits consisting in rotating waves $u(x,t)=u(x-ct)$
 (notice the
correspondence with the closed orbits of a radially
symmetric two-dimensional flow).
This property is a consequence of the
zero number property of Theorem
\ref{th-Angenent} and has been proved in \cite{Angenent-Fiedler} by 
Angenent and
Fiedler. The works of Matano and Nakamura \cite{Matano-Nakamura} and of
Fiedler, Rocha and Wolfrum \cite{Fiedler-Rocha-Wolfrum} show that the unstable
and stable manifolds of the equilibria and the periodic orbits always
intersect transversally and that no homoclinic orbit can occur. Moreover,
in \cite{Fiedler-Rocha-Wolfrum}, the authors give
an algorithm for determining the global dynamics of the parabolic equation
\eqref{PDE} on $\Omega=\Tm^1$ with homogeneous nonlinearity $f\equiv
f(u,\partial_x u)$. This algorithm uses the knowledge of the equilibria and the
periodic orbits only.
In \cite{Rocha}, Rocha also characterized all the dynamics, which may
occur. Due to the automatic 
transverse intersection of the stable and unstable manifolds
and due to the
possibility of transforming any circle of equilibrium points into a rotating
periodic orbit (see
\cite{Fiedler-Rocha-Wolfrum}), one can show that the Morse-Smale property
holds for the parabolic equation on $\Tm^1$ for a generic homogeneous
nonlinearity $f(u,\partial_x u)$
(see \cite{JR}). Finally, the realization of any radially symmetric
two-dimensional flow in the dynamics of \eqref{PDE} on $\Tm^1$
for some $f\equiv f(u,\partial_x u)$ and the fact that the dimension of the
attractor is equal to the largest dimension of the unstable manifolds are shown
in \cite{HJR}.

\vspace{2mm}

{\noindent \bf $\bullet$ $d\geq 3$ and dim($\Omega$)$\geq 2$}\\
The genericity of the Kupka-Smale property for vector fields in $\Rm^d$, $d\geq 3$, has
been proved independently by Kupka in \cite{Kupka} and by Smale in
\cite{Smale2}. Their proofs have been simplified by Peixoto in
\cite{Peixoto2} (see \cite{Abraham-Robin} and \cite{Palis-de-Melo}). The strong
difference with the lower dimensional vector fields is
that, when $d\geq 3$, \eqref{ODE} may admit transversal homoclinic orbits
consisting in the transversal intersection of the stable and unstable manifolds
of a hyperbolic periodic orbit. The existence of such an intersection is stable
under small perturbations and yields a Smale horseshoe containing an infinite
number of
periodic orbits and chaotic dynamics equivalent to the dynamics of the shift
operator, see \cite{Smale65}. Therefore, the Morse-Smale property cannot be
dense in vector fields on $\Rm^d$ with $d\geq 3$.
Even worse, the set of vector fields, whose dynamics are structurally
stable under small perturbations, is not dense (notice that this set contains
the vector fields satisfying the Morse-Smale property). Indeed, as shown in
\cite{Guckenheimer-Williams}, there exists an open set $\Um$ of vector fields
of $\Rm^3$ and a foliation of $\Um$ by
$2-$codimensional leaves $(\Um_\lambda)_{\lambda\in\Rm^2}$ such that
each $g\in\Um$ admits a Lorenz attractor $\Am_g$ and such that the dynamics of
two attractors $\Am_g$ and $\Am_{\tilde g}$ are qualitatively
equivalent if and only if $g$ and $\tilde g$ belong to the same leave
$\Um_\lambda$.
The possible presence of other chaotic dynamics such as Anosov systems or wild
dynamics is also noteworthy, see \cite{Anosov}, \cite{Newhouse} and
\cite{Bonatti-Diaz}. For the interested reader, we refer to \cite{Robinson2} or
\cite{Smale67}.

In \cite{BJR}, Brunovsk\'y and both authors proved that the stable and 
unstable manifolds of hyperbolic equilibria or periodic orbits of the 
parabolic equation \eqref{PDE} are generically 
transversal. To obtain the genericity of Kupka-Smale property, it remains to 
obtain the generic hyperbolicity of periodic orbits. Even if this problem is 
still open, we may strongly believe that the Kupka-Smale property is generic 
for the parabolic equation \eqref{PDE}.
There exist several results concerning the
embedding of the finite-dimensional flows into the parabolic equations.
Pol\'a\v{c}ik has shown in \cite{Polacik5} that any ordinary differential
equation on $\Rm^d$ can be embedded into the flow of \eqref{PDE} for some $f$
and for some domain $\Omega\subset\Rm^d$. The constraint that the dimension of
$\Omega$ is equal to the dimension of the imbedded flow is removed in
\cite{Polacik}, however the result concerns a dense set of flows 
only. A similar
  result has been obtained by Dancer and Pol\'a\v{c}ik in
\cite{Dancer-Polacik} for homogeneous nonlinearities $f(u,\grad u)$ (see also
\cite{Prizzi-Rybakowski}). These realization results imply the possible
existence of
persistent chaotic dynamics in the flow of the parabolic equation 
\eqref{PDE} as
soon as
$\Omega$ has a dimension larger than one: transversal homoclinic orbits, Anosov
flows on invariant manifolds of any dimension, Lorenz attractors etc...

\vspace{2mm}

{\noindent \bf $\bullet$ Gradient case}\\
When $g$ is a gradient vector field $\grad G$ with $G\in\Cm^2(\Rm^d,\Rm)$, then
$-G$ is a strict Lyapounov function and \eqref{ODE} is a gradient system.
In this case, the Kupka-Smale property is equivalent to the Morse-Smale
property. The
genericity of the Morse-Smale property for gradient vector fields has been obtained
by Smale in \cite{Smale1}.

In the case where the nonlinearity $f \in\Cm^r(\overline\Omega \times\Rm, \Rm)$ (that is,
 $f\equiv f(x,u)$ does not depend on $\grad u$),
the parabolic equation \eqref{PDE} admits a strict Lyapounov function
 given by $E(u)=\int_\Omega\left(\frac 12 |\grad u(x)|^2 - F(x,u(x))\right) dx$,
where
$F(x,u)=\int_0^u f(x,s)ds$ is a primitive of $f$, and hence generates
a gradient system. Brunovsk\'y and Pol\'a\v{c}ik have shown in 
\cite{Bruno-Pola} that the
Morse-Smale property holds for the parabolic equation,
generically with respect to $f(x,u)$. It is noteworthy that
the Morse-Smale property is no longer generic if one restricts  the
nonlinearities  
to the class of homogeneous functions $f\equiv f(u)$  (see
\cite{Polacik2}). Pol\'a\v{c}ik has shown in \cite{Polacik3} that any
generic gradient vector field of $\Rm^d$ can be realized in the flow of the
parabolic equation \eqref{PDE} on a bounded domain of $\Rm^2$ with an appropriate nonlinearity
$f(x,u)$. The paper \cite{Polacik3} also contains the realization of particular dynamics such
as non-transversal intersections of stable and unstable manifolds.

\vspace{2mm}

{\noindent \bf $\bullet$ Caveat: general ODEs or cooperative systems?}\\
In Table 1, we have given the striking correspondence
between the flow
generated by Eq. \eqref{ODE} and the semiflow generated
by the parabolic equation
\eqref{PDE}. In addition, we have pointed out that some classes of cooperative systems are
also involved in this correspondence. In fact, the reader should be aware that
the dynamics of \eqref{PDE} is much closer to the ones of  a cooperative
system than to the ones of the general vector field 
\eqref{ODE}. Indeed, the semiflow
$S_f(t)$ generated by the parabolic equation \eqref{PDE} belongs to the class of
strongly monotone semiflows, which means that this semiflow has more constraints than the flow
$T_g(t)$ generated by a general vector field $g$ (see Section \ref{section-coop}).
That is why, it could be more relevant to write Table 1  in terms of
cooperative systems  only (for example, by replacing the case of the
general ODE with $d\geq 3$ by the case of a cooperative system of ODEs in
dimension $d\geq 4$). However, we have chosen to mainly write Table 1 in terms
of general ODEs for several reasons:\\
- as far as the properties stated in Table 1 are concerned, there is no
difference between the dynamics of a general ODE and the ones of a parabolic
PDE,\\
- the dynamics of general ODEs are common knowledge, whereas speaking in terms
of cooperative systems may not give a good insight of the dynamics of
\eqref{PDE},\\
- not all the properties stated in Table 1 are known for the class of
cooperative systems (for example the genericity of
Kupka-Smale property is not yet known for $d\geq 4$).

\section{Zero number and unique continuation properties for the 
scalar parabolic
equation}\label{section-ucp}

The results presented in Table \ref{table} and in Section
\ref{section-comments} strongly rely on properties specific to the
parabolic equations.
The purpose of this section is to give a first insight of these
particular properties and of their use, to the reader.

Dynamical systems generated by vectors fields are flows on $\Rm^d$, whereas the
phase-space of the
  parabolic equation is an infinite-dimensional space
$X^\alpha$.  It is important to be aware of the fact that the
parabolic equations generate only a small part of all possible
dynamical systems on the Banach space $X^\alpha$. On one hand, this 
implies less
freedom in
perturbing the dynamics and hence in obtaining density results. In particular,
whereas one can easily construct perturbations of a vector field $g$ which
are localized in the phase space $\Rm^d$, the perturbations of the nonlinearity
$f$ act in a non local way on $X^\alpha$ (many different functions
$u$ can have the same values of $u$ and $\grad u$ at a given point $x$).
Therefore, it is important to obtain unique continuation results in order to
find values $(x,u,\grad u)$, which are reached only once by a given periodic,
heteroclinic or homoclinic orbit. On the other hand, the small class of
dynamics generated by the parabolic equations admits special properties. These
properties may in particular yield the constraints, which make the dynamics
similar to the ones of low-dimensional vector fields.

\vspace{2mm}

The scalar parabolic equation in space dimension one ($\Omega=(0,1)$
or $\Tm^1$) satisfies a very strong property: the number of zeros of the
solutions of the linearized equation is nonincreasing in time. This property is
often called Sturm property since its idea goes back to Sturm \cite{Sturm} in
1836. There are different versions of this result, which have been proved by
Nickel \cite{Nickel}, Matano \cite{Matano,Matano2}, Angenent and Fiedler
\cite{Angenent,Angenent-Fiedler} and Chen \cite{Chen} (see also
\cite{GalakHarwin} for
a survey).
By similar technics, a geometrical result on braids
formed by solutions of the one-dimensional parabolic equation is obtained in
\cite{GV}.
\begin{theorem}\label{th-Angenent}
Let $\Omega=(0,1)$ with Neumann boundary conditions or $\Omega=\Tm^1$.
Let $T>0$, $a\in W^{1,\infty}(\overline\Omega\times [0,T],\Rm)$ and
$b\in\Lm^\infty(\overline\Omega\times [0,T] ,\Rm)$.
Let $v:\overline\Omega\times (0,T)\rightarrow \Rm$ be a bounded
non-trivial classical solution of
$$\partial_t v=\partial^2_{xx}v+a(x,t)\partial_x
v+b(x,t)v~~,~~(x,t)\in\Omega\times (0,T)~.$$
Then, for any $t\in (0,T)$, the number of zeros of
the function $x\in\overline\Omega\mapsto
v(x,t)$ is finite and non-increasing in time. Moreover, it strictly
decreases at $t=t_0$ if and only if $x\mapsto v(x,t_0)$ has a
multiple zero.
\end{theorem}
Theorem \ref{th-Angenent} is the fundamental ingredient of almost all the
results given in Table \ref{table} in the cases
$\Omega=(0,1)$ and $\Omega=\Tm^1$. It can be used either as a strong comparison
principle or as a strong unique continuation property, as shown in the
following examples of applications. General surveys can be found in
\cite{Fiedler-Scheel}, \cite{HaleCanada} and \cite{HJR}.

In the first application presented here, Theorem \ref{th-Angenent} is used as a
strong
maximum principle. In some sense, it yields an order on the
phase space which is preserved by the flow. This illustrates how
Theorem \ref{th-Angenent} may imply constraints similar to the ones of
low-dimensional vector fields. The following result was first proved in
\cite{Zelenyak} and the proof given here comes from \cite{Matano} (see also
\cite{Fiedler-Scheel}).
\begin{prop}\label{prop-cv-single}
Let $\Omega=(0,1)$, let $u_0\in X^\alpha$ and let $u(x,t)$ be the
corresponding solution of the parabolic equation \eqref{PDE} with homogeneous
Neumann boundary conditions. The $\omega-$limit
set of $u_0$ consists of a
single
equilibrium point.
\end{prop}
\begin{demo}
We first notice that $v(x,t)=\partial_t u(x,t)$ satisfies the equation
$$\partial_t v(x,t)=\partial^2_{xx} v(x,t)+ f'_{u}(x,u(x,t),\partial_x u(x,t))
  v(x,t) + f'_{\partial_x u}(x,u(x,t),\partial_x u(x,t)) \partial_x v(x,t)~.$$
Due to the Neumann boundary conditions, we have $\partial_x
u(0,t)=\partial^2_{xt} u(0,t)=\partial_x v(0,t)=0$ for all
$t >0$.
 In particular, as soon as $v(0,t)=0$, $v(t)$ has a double zero at
$x=0$. Due to Theorem \ref{th-Angenent}, either $v$ is a trivial solution, that
is $v\equiv 0$ for all $t$, and $u$ is an equilibrium
point, or $v(0,t)$ vanishes at most a finite number of times since $v(t)$
can have a multiple zero only a finite number of times.
Assume that $u$ is not an equilibrium, then $u(0,t)$ must be monotone for
large times and thus converges to $a\in\Rm$.
Any trajectory $w$ in
the $\omega-$limit set of $u_0$ must hence satisfy $w(0,t)=a$ for all
$t$. Therefore, $\partial_t w(0,t)=0$ for all $t$ and $\partial_t
w(0,t)$ has a multiple zero at $x=0$ for all times. Using Theorem
\ref{th-Angenent}, we deduce as above that $w$ is an equilibrium
point of \eqref{PDE}. But there exists at most one equilibrium $w$
satisfying $w(0)=a$ and the Neumann boundary condition $\partial_x w(0)=0$.
Therefore, the $\omega-$limit set of $u_0$ is a single equilibrium point $w$.
\end{demo}

The second application comes from \cite{JR}. It shows how Theorem
\ref{th-Angenent} can be used as a unique continuation property. This kind of
property roughly says that, if two solutions coincide too much
near a point $(x_0,t_0)$, then they must
be equal everywhere.
The motivation beyond this example of application is the following.
We consider a time-periodic
solution of \eqref{PDE} on $\Omega=\Tm^1$. The problem is to find a 
perturbation
of the nonlinearity $f$, which makes this periodic orbit hyperbolic. 
As enhanced
above, such a perturbation is nonlocal in the phase space of \eqref{PDE}. To be
able to perform perturbation arguments, it is important to show that one can
find a perturbation of $f$ which acts only locally on the periodic orbit. To
this end, one proves the following result.
\begin{prop}\label{prop-injectivite}
Let $p(x,t)$ be a periodic orbit of \eqref{PDE} on $\Omega=\Tm^1$. Let $T>0$ be
its minimal period. Then, the map
$$(x,t) \in \Tm^1 \times [0, T) \mapsto (x, p(x,t), \partial_x p(x,t))
$$
is one to one.
\end{prop}
\begin{demo}
Assume that this map is not injective. Then there exist
$x_0$, $t_0 \in [0, T)$ and $t_1 \in [0, T)$, $t_0\ne t_1$ such
that
$$ p(x_0, t_0)=p(x_0, t_1)~\text{ and }~\partial_x p(x_0, t_0)=
\partial_x p(x_0, t_1)~.$$
The function $v(x,t)= p(x, t + t_1 -t_0)-p(x,t)$ is a solution of
the equation
$$
\partial_t v(x,t)= \partial^2_{xx} v(x,t) + a(x,t)v(x,t) + b(x,t)
\partial_x v(x,t)~,
$$
where $a(x,t)= \int_0^1 f'_u(x,p(x,t) + s(p(x,t+
t_1-t_0)-p(x,t)), \partial_x p(x,t + t_1-t_0)) ds$ and
$b(x,t) =\int_{0}^{1} f'_{u_x}(x,p(x,t), \partial_x (p(x,t) +s(p(x,t
+ t_1-t_0)-p(x,t))) ds$. Moreover, the function $v(x,t)$ satisfies
$v(x_0,t_0)=0$ and $\partial_x v(x_0, t_0)=0$ and does not vanish
everywhere since $|t_1-t_0|<T$. Due to Theorem \ref{th-Angenent}, the
number of zeros of $v(t)$ drops strictly at $t=t_0$ and never increases.
However, $v(t)$ is a periodic function of period $T$, and thus, its number
of zeros is periodic. This leads to a contradiction and proves the proposition.
\end{demo}

\vspace{2mm}

In a domain $\Omega$ of dimension $d'\geq 2$, there is no known
counterpart for Theorem \ref{th-Angenent} as shown in \cite{Fusco-Lunel}. In
particular, Proposition \ref{prop-injectivite} does no longer hold. However, to
be able to construct relevant perturbations of periodic orbits, one needs a
result similar to Proposition \ref{prop-injectivite}, even if weaker. The
following result can be found in \cite{BJR}. Its proof is based on a
generalization of the arguments of \cite{Hardt-Simon} and on unique
continuations properties of the parabolic equations.
\begin{theorem}\label{th-nodal-sing}
Let $p(x,t)$ be a periodic orbit of \eqref{PDE} with minimal period $T>0$.
There exists a generic set of points
$(x_0,t_0)\in \Omega \times [0,T)$ such that if $t\in [0,T)$ satisfies
$p(x_0,t)=p(x_0,t_0)$ and $\grad p(x_0,t)=\grad p(x_0,t_0)$, then $t=t_0$.
\end{theorem}


\section{Cooperative systems of ODEs}\label{section-coop}
We consider a system of differential equations
\begin{equation}\label{eq-coop}
\dot y(t)=g(y(t))~, \quad y(0)=y_0\in\Rm^N~,
\end{equation}
where $g=(g_i)_{i=1\ldots N}$ is a $\Cm^1$ vector field.\\
Due to the analogy with biological models, the following definitions are
natural. We say that \eqref{eq-coop} is a {\it cooperative} (resp.
{\it competitive}) system if for any $y\in\Rm^N$ and $i\neq j$, $\Drond
{g_i}{y_j}(y)$ is non-negative (resp. non-positive) 
and the matrix
$(\Drond{g_i}{y_j})(y)$ is irreducible i.e. it is not a block diagonal matrix
(the simpler assumption that all the coefficients $\Drond {g_i}{y_j}(y)$ are
positive is sometimes made instead of the irreducibility).
We say that \eqref{eq-coop} is a {\it tridiagonal} system if
$\Drond {g_i}{y_j}=0$ for $|i-j|\geq 2$ and a {\it cyclic tridiagonal} system
if the indices $i$ and $j$ are considered modulo $N$, i.e. if, in addition, we
allow $\Drond {g_1}{y_N}$ and $\Drond {g_N}{y_1}$ to be non-zero. For
the reader interested in cooperative systems, we refer to \cite{Smith2}.

In this section, we only consider cooperative systems. However, notice that, by
changing $t$ into $-t$ or $y_i$ into $-y_i$, we obtain similar results for
competitive systems and for systems with different sign
conditions.

The dynamics of cooperative systems may be as complicated as the
dynamics of general vector fields. Indeed, Smale has shown in \cite{Smale-coop}
that any vector field in $\Rm^{N-1}$ can be realized in a invariant
manifold of a cooperative system in $\Rm^N$. 
Notice that this realization
result implies that any one-dimensional vector field can be imbedded in a
tridiagonal cooperative system and any two-dimensional vector field can be
imbedded in a cyclic tridiagonal cooperative system. This explains why we
present the tridiagonal cooperative systems in Table \ref{table} as
generalization of one- and two-dimensional vector fields.

However, the dynamics of a cooperative system \eqref{eq-coop} is really
different from the ones of the general ODE \eqref{ODE} since a cooperative system
generates a strongly monotone flow, that is, a flow which preserves a partial
order. It is noteworthy that the semiflow $S_f(t)$ generated by the parabolic
equation \eqref{PDE} also belongs to the class of strongly
monotone semiflows (it preserves the order of $X^\alpha$ induced by the
classical order of $\Cm^0(\Omega)$). Therefore, the semiflow of \eqref{PDE} is
much closer to the one of  the cooperative system \eqref{eq-coop}, both admitting
more constraints than the flow $T_g(t)$ generated by a general vector field
$g$. In \cite{Hirsch0} and \cite{Hirsch3} for
example, Hirsch has shown  that almost all bounded trajectories of a strongly
monotone semiflow are quasiconvergent, that is, their $\omega$-limit
sets consist only of equilibria. More precisely, all initial data, which
have bounded nonquasiconvergent trajectories, form a meager subset 
(that is, the  complement of a generic subset) of the phase space. Later,
in \cite{Polacik0}, Pol\'a\v{c}ik has proved that the set of all initial data
$u_0 \in X^{\alpha}$, which have bounded nonconvergent trajectories in the
semiflow of the parabolic equation \eqref{PDE}, is meager in $X^{\alpha}$.

Moreover, since the works of Hirsch
and Smillie, it is known that the dynamics of cooperative systems, which are in
addition tridiagonal, are very constrained in any dimension $N$. Indeed, in
\cite{Hirsch1},
\cite{Hirsch2} and \cite{Smillie}, strong properties of the limit sets of
cooperative systems are proved. In particular, any three-dimensional 
cooperative
system satisfies the Poincar\'e-Bendixson property and the trajectory of any
tridiagonal cooperative system converges to a single equilibrium 
point. Inspired
by the articles of Henry and Angenent about the parabolic equation on $(0,1)$,
Fusco and Oliva (see \cite{Fusco-Oliva1}) showed a theorem similar to Theorem
\ref{th-Angenent} (see \cite{Smith} for a more general statement).

\begin{theorem}\label{th-FO}
Let $\Nc$ be the set of vector $y\in\Rm^N$ such that, for all $i=1\ldots N$, either $y_i\neq 0$ or $y_i=0$ and $y_{i-1}y_{i+1}<0$ (where $y_0=y_{N+1}=0$).
For every $y\in\Nc$, we set $N(y)$ to be the number of sign changes
for $y_i$,  when $i$ goes from $1$ to $N$.
Let $y(t)\neq 0$ be a solution of
\begin{equation}\label{eq-th-FO}
\dot y(t)=A(t)y(t)~,
\end{equation}
where $A\in\Cm^0(\Rm,\Mc_N(\Rm))$ satisfies $A_{ij}(t)>0$ for all $t\in\Rm$ and
all $i\neq j$.\\
Then, the times $t$ where $y(t)\not\in \Nc$ are isolated and, if
$y(t_0)\not\in\Nc$, then, for every $\varepsilon>0$ small enough,
$N(y(t+\varepsilon))<N(y(t-\varepsilon))$.
\end{theorem}
In other words, the number of sign changes of the solutions of the linear
equation \eqref{eq-th-FO} is non-increasing in time and strictly drops at
$t_0$ if and only if $y(t_0)$ has in some sense a multiple zero. The parallel
with Theorem \ref{th-Angenent} is of course striking. Using Theorem
\ref{th-FO}, Fusco and Oliva have shown that the stable and unstable manifolds
of equilibrium points of a tridiagonal cooperative system always intersect
transversally. As a consequence, the Morse-Smale property is generic in the
class of tridiagonal cooperative systems.

Theorem \ref{th-FO} also holds for cyclic tridiagonal cooperative systems, see
\cite{Fusco-Oliva2} and \cite{Smith}. Using this fundamental property, Fusco
and Oliva have shown in \cite{Fusco-Oliva2} that  the stable and unstable
manifolds of periodic orbits of cyclic tridiagonal cooperative systems always
intersect transversally.
In addition, Mallet-Paret and Smith have
shown in \cite{MP-Smith} that cyclic tridiagonal cooperative systems satisfy
the Poincar\'e-Bendixson property. Notice that, following \cite{JR} and
\cite{JR2}, one should be able to prove the genericity of the Morse-Smale
property for cyclic tridiagonal cooperative systems. 
This has been proved very recently by Percie du Sert (see \cite{Percie})\\[3mm]

Considering all these results, it is not surprising that there exists a
parallel between tridiagonal cooperative systems and the parabolic equation on
$(0,1)$. Indeed, consider a solution $v$ of the linear
one-dimensional parabolic equation
\begin{equation}\label{eq-lien}
\dot v(x,t)=\partial^2_{xx} v(x,t) + a(x,t)\partial_x v(x,t) +
b(x,t)v(x,t)~~~ (x,t)\in\Rm\times\Rm_+
\end{equation}
We discretize the segment $(0,1)$ by a sequence of points 
$x_k=(k-1)/(N-1)$ with
$k=1\ldots N$. The natural approximation of $v$ is given by $y_k\approx
v(x_k)$ solution of
\begin{equation}\label{eq-lien2}
\dot y_k(t)=\frac {y_{k+1}(t)-2y_k(t)+y_{k-1}(t)}{h^2} +
a_k(t)\frac{y_{k+1}-y_k}h  + b_k(t)y_k(t)~~~
\end{equation}
where $a_k(t)=a(x_k,t)$, $b_k(t)=b(x_k,t)$ and $h=1/(N-1)$. If $h$ is small
enough, \eqref{eq-lien2} is a tridiagonal cooperative system. The relation
between Theorems \ref{th-Angenent} and \ref{th-FO} is obvious in this 
framework.


\section{Kupka-Smale property and other open problems}\label{section-open}

One of the main goals of the study of dynamical systems is to understand the
behaviour of a generic dynamical system. The most recent results
concerning the parabolic equations are related to the genericity of Kupka-Smale
property. In \cite{BJR}, the generic transversality of homoclinic and 
heteroclinic orbits of hyperbolic equilibria or periodic orbits of \eqref{PDE} 
is proved. However, the generic hyperbolicity of 
periodic orbits is still an open problem (see the discussion in \cite{BJR}). 
Its resolution would complete the whole correspondance of Table 1 and is one of 
the most important open problem concerning the qualitative dynamics of the 
parabolic equation \eqref{PDE}.

However, even if the generic Kupka-Smale property is obtained, it 
cannot give a good insight of the complex and chaotic dynamics that may be 
generated by homoclinic connections. For this reason, the study of 
finite-dimensional flows has been pursued further the Kupka-Smale property and 
is still in progress. The corresponding results should serve as a guideline for
the study of the flow generated by the parabolic equation \eqref{PDE}.
For vector fields, one of the main steps beyond Kupka-Smale property is Pugh's 
closing lemma: if $p$ is a non-wandering point
of the dynamical system $T_g(t)$ generated by \eqref{ODE}, then there exists a
$\Cm^1-$perturbation $\tilde g$ of $g$ such that $p$ is a periodic
point of $T_{\tilde g}(t)$ (the case of a $\Cm^r-$perturbation with $r\geq2$
is still open). The proof of Pugh in \cite{Pugh} concerns discrete dynamical
systems. It has been adapted to the case of flows by Pugh and Robinson in
\cite{Pugh-Robinson} (see also \cite{Robinson3} for an introduction to the
proof). A direct consequence of Pugh closing lemma is the general density
theorem: for a generic finite-dimensional flow, the non-wandering points are is
the closure of the periodic points (see \cite{Pugh2} and \cite{Robinson2}).
Other connecting lemmas have been proved by Hayashi \cite{Hayashi} and Bonatti
and Crovisier \cite{Bonatti-Crovisier}. They enable a better understanding of
generic dynamics. For example, the class of finite-dimensional 
dynamical systems
which either satisfy the Morse-Smale property or admit a transversal homoclinic
connection is generic (see \cite{Pujals-Sambarino}, \cite{Bonatti-Gan-Wen} and
\cite{Crovisier} for discrete dynamical systems in dimensions $d=2$, $d=3$
and $d\geq 4$ respectively, and see \cite{Arroyo-Rodriguez} for
three-dimensional flows).
Obtaining similar results for the flow of the parabolic equation should be a
very interesting and difficult challenge.

Other interesting open problems concern the realization of finite-dimensional
dynamics in the semiflow of parabolic equations. Indeed, we only know
that one can realize the dynamics of a dense set of general ODEs in
the flow of a parabolic equation \eqref{PDE} on a two-dimensional domain. One
may wonder if it is possible to realize the dynamics of all ODEs. Since the
parallel between parabolic equations and cooperative systems is stronger,
the following strong realization conjecture may be more plausible: \emph{any}
flow of a cooperative system of ODEs can be realized in an invariant
manifold of the flow of a parabolic equation \eqref{PDE} on a two-dimensional
domain.

Finally, the genericity of  the Morse- and Kupka-Smale properties 
is also an interesting problem for other classes of partial differential
equations. The genericity of the Morse-Smale property is known for the wave
equations $\ddot u+\gamma \dot u =\Delta u+f(x,u)$ with constant damping
$\gamma>0$ (see \cite{Bruno-Raugel}) and with variable damping 
$\gamma(x)\geq 0$
in space dimension one (see \cite{Joly}).
We recall that, in both cases, the associated dynamical system is gradient. 
Nothing is known for other classes of
PDEs, in particular for the equations of fluids dynamics and for systems of
parabolic
equations $\dot U=\Delta U+f(x,U)$, with $U(x,t)\in\Rm^N$.
In all these cases, the main problem consists in understanding how the perturbations
act on the phase plane of the PDE. Either one proves unique
continuation results similar to Theorem \ref{th-nodal-sing} in order to be able to use
local perturbations of the flow (as in \cite{JR}, \cite{JR2} and 
\cite{BJR}), or
one uses particular non-local perturbations in a very careful way (as in
\cite{Bruno-Pola}, \cite{Bruno-Raugel} and \cite{Joly}).


\section*{Glossary}
In this section, $S(t)$ denotes a general continuous dynamical system
on a Banach space $X$.
An orbit of $S(t)$ is denoted by $x(t)=S(t)x_0$ with $t\in I$, where
$I=[0,+\infty)$, $I=(-\infty,0]$ or $I=(-\infty,+\infty)$ in the case of a
positive, negative or global trajectory respectively.

{\noindent\bf Compact global attractor:} if it exists, the compact
global attractor
$\Am$ of $S(t)$ is a compact invariant set which attracts all the bounded sets
of $X$. Notice that $\Am$ is then the set of all the bounded global
trajectories. See \cite{Hal88}.

{\noindent\bf $\alpha-$ and $\omega-$limit sets:} let $x_0\in X$.
The $\alpha-$limit set $\alpha(x_0)$ and the $\omega-$limit set
$\omega(x_0)$ of $x_0$  are the sets of accumulation points of
  the negative and positive orbits coming from $x_0$
respectively. More precisely,
\begin{align*}
\alpha(x_0)&=\{ x\in
X~/~\exists (t_n)_{n\in\Nm},~t_n\xrightarrow[n\rightarrow
\infty]{}-\infty~\text{
   and a negative trajectory
} x(t)\\ &~~~~\hbox{such that }x(0)=x_0~\hbox{ and }~x(t_n)
\xrightarrow[n\rightarrow
\infty]{} x \} \\
\omega(x_0)&=\{ x\in X~/~\exists
(t_n)_{n\in\Nm},~t_n\xrightarrow[n\rightarrow \infty]{}+\infty~\text{
   such that }S(t_n)x_0\xrightarrow[n\rightarrow \infty]{} x \}
\end{align*}
The limit sets $\alpha(x_0)$ and $\omega(x_0)$ are non-empty connected compact
sets.

{\noindent\bf Homoclinic or heteroclinic orbit:} let $x(t)=S(t)x_0$ be a
global trajectory of $S(t)$. Assume that the $\alpha-$ and
$\omega-$limit sets of $x_0$ exactly consists in one orbit,
denoted $x_-(t)$ and $x_+(t)$ respectively, this
orbit being either an equilibrium point or a periodic orbit. The
trajectory $x(t)$
is said to be a homoclinic orbit if $x_-(t)=x_+(t)$ and a heteroclinic
orbit if $x_-(t)\neq x_+(t)$.

{\noindent\bf Backward uniqueness property:} $S(t)$ satisfies the backward
uniqueness property if for any time $t_0>0$ and any trajectories $x_1(t)$ and
$x_2(t)$, $x_1(t_0)=x_2(t_0)$ implies $x_1(t)=x_2(t)$ for all $t\in [0,t_0]$.
Notice that this does not mean that $S(t)$ admits negative trajectories.

{\noindent\bf Hyperbolic equilibrium points or periodic orbits:} an equilibrium
point $e$ of $S(t)$ is hyperbolic if the linearized operator 
$x\mapsto D_eS(1)x$
has no spectrum on the unit circle. 
Let $p(t)$ be a  periodic solution of $S(t)$
with minimal period $T$. For each $x\in X$, we denote $t\mapsto\Pi(t,0)x$ the
corresponding  trajectory of the linearization of $S(t)$ along the periodic
solution $p(t)$. Then, $p(t)$ is said hyperbolic if the linear map
$x\mapsto \Pi(T,0)x$ has no spectrum on the unit circle except the eigenvalue
$1$ which is simple. Remark  that then, for any integer $k \ne 0$, the linear 
map $x\mapsto \Pi(kT,0)x$ has no spectrum on the unit circle except the eigenvalue
$1$ (which is simple).

{\noindent\bf Stable and unstable manifolds:} let $e$ be a hyperbolic
equilibrium point of $S(t)$. There exists a neighbourhood $\Nc$ of 
$e$ such that
the set
\begin{equation*}
\begin{split}
W^u_{loc}(e)=\{x_0\in X~,~\exists\hbox{ a negative 
trajectory }x(t)\hbox{ with }x(0&)=x_0 \cr
&\hbox{ and, } \forall t\leq 0,~x(t)\in\Nc\}
\end{split}
\end{equation*}
is a
submanifold of $X$,
in which all negative trajectories converge to $e$ when $t$
goes to $-\infty$. The manifold $W^u_{loc}(e)$ is called the local unstable
manifold of $e$. Pushing $W^u_{loc}(e)$ by the flow $S(t)$, one  can define the
(global) unstable set $W^u(e)=\cup_{t\geq 0} S(t)W^u_{loc}(e)$, which consists
in all the negative trajectories converging to $e$ when $t$ goes to $-\infty$.
This unstable set $W^u(e)$ is
an immersed submanifold
under suitable properties. For
instance, backward uniqueness properties are needed to extend the manifold
structure.\\
In the same way, one defines the local stable manifold
\begin{align*}
W^s_{loc}(e)&=\{x_0\in X~,~\forall t\geq 0,~S(t)x_0\in\Nc\}\\
&=\{x_0\in X~,~\forall t\geq 0,~S(t)x_0\in\Nc\text{ and
}S(t)x_0\xrightarrow[t\rightarrow +\infty]{}e\}~.
\end{align*}
General partial differential equations (and parabolic equations in 
particular) do not admit negative trajectories.
Therefore, it is less easy to extend the  local 
stable manifold to a global stable manifold. However, one can define the stable set 
$W^s(e)$ of $e$ as follows
$$
W^s(e) =\{x_0\in X~, ~S(t)x_0\xrightarrow[t\rightarrow +\infty]{}e\}~.
$$
Under suitable additional properties (which are satisfied by the parabolic equation \eqref{PDE}), 
one can show that $W^s(e)$ is also an immersed submanifold. For
instance, backward uniqueness properties of the adjoint dynamical system $S^*(t)$ on $X^*$ and 
finite-codimensionality of $W^s_{loc}(e)$ are needed (see \cite{Henry-book} for more details).\\
If $p(t)$ is a hyperbolic periodic orbit, one defines its unstable and local
stable manifolds in the same way.
See for example \cite{Palis-de-Melo} for more details.

{\noindent\bf Non-wandering set:} a point $x_0\in X$ is non-wandering 
if for any
neighbourhood $\Nc\ni x_0$ and any time $t_0>0$, there exists $t\geq t_0$ such
that $S(t)\Nc\cap \Nc\neq\emptyset$.

{\noindent\bf The Kupka-Smale and Morse-Smale properties:} $S(t)$ satisfies the
 Kukpa-Smale property if all its equilibrium points and  periodic orbits are
hyperbolic and if their stable and unstable manifolds intersect transversally.
It satisfies the Morse-Smale property if in addition its non-wandering set
consists only in a finite number of equilibrium points and periodic 
orbits. We refer
to \cite{Palis-de-Melo} for more precisions on these notions.

{\noindent\bf Gradient dynamical systems:} $S(t)$ is gradient if it admits a
Lyapounov functional, that is a function $\Phi\in\Cm^0(X,\Rm)$ such 
that, for all
$x_0\in X$, $t\mapsto \Phi(S(t)x_0)$ is non-increasing and is constant if and
only if $x_0$ is an equilibrium point. We recall that a gradient dynamical
system does not admit periodic or homoclinic orbits.

{\noindent\bf Cooperative system of ODEs:} see Section \ref{section-coop}.

{\noindent\bf Generic set and Baire space:} a generic subset of a
topological space $Y$ is a set which contains a countable intersection of dense
open subsets of $Y$. A property is generic in $Y$ if it is satisfied for a
generic set of $Y$. The space $Y$ is called a Baire space if any generic set is
dense in $Y$. In particular a complete metric space is a Baire space.

{\noindent\bf Whitney topology:} let $k\geq 0$ and let $M$ be a 
Banach manifold.
The Whitney topology on $\Cm^k( M,\Rm)$ is the topology generated by the
neighbourhoods
$$ \{g\in\Cm^k( M,\Rm)~,~|D^i f(x)-D^i g(x)|\leq\delta(x),~\forall
i\in\{0,1,\ldots,k\},~
\forall x\in M \}~,$$
where $f$ is any function in $\Cm^k( M,\Rm)$ and $\delta$ is any positive
function in $\Cm^k( M,(0,+\infty))$.
Notice that $\Cm^k( M,\Rm)$ endowed with the Whitney topology is a 
Baire space even
if it is not a metric space when $M$ is not compact. We refer for instance to
\cite{GG}.

{\noindent\bf The fractional power space $X^\alpha$:} let $A$ be a positive
self-adjoint operator with compact inverse on $\Lm^2(\Omega)$. Let 
$(\lambda_n)$
be the sequence of its eigenvalues, which are positive, and let
$(\varphi_n)$ be
the corresponding sequence of eigenfunctions, which is an orthonormal basis of
$\Lm^2(\Omega)$. For each $\alpha\in\Rm$, we define the fractional power of $A$
by $A(\sum_n c_n \varphi_n)=\sum_n c_n \lambda_n^\alpha \varphi_n$. In
particular, $A^0=Id$ and $A^1=A$. The space $X^\alpha$ is the domain of
$A^\alpha$ that is
$$X^\alpha=\{\varphi\in\Lm^2(\Omega)~,~\varphi=\sum_n c_n\varphi_n \hbox{ such
that }(c_n\lambda_n^\alpha)\in\ell^2(\Nm) \}~.$$
It is possible to define fractional powers of more general operators, called
sectorial operators, see \cite{Henry-book}.

{\noindent\bf The Sobolev space $W^{s,p}(\Omega)$:} if $s$ is a positive
integer, $W^{s,p}(\Omega)$ is the space of (classes of) functions $f\in\Lm^p(\Omega)$, which
are $s$ times differentiable in the sense of distributions and whose derivatives up to order
 $s$ belong to $\Lm^p(\Omega)$. It is possible to extend
this notion to positive numbers $s$ which are not integers by using 
interpolation theory.

{\noindent\bf Unique continuation properties:} let us consider a partial
differential equation on $\Omega$ and let $u(x,t)$ be any solution
of it. A
unique continuation property for this PDE is a result stating that if $u(x,t)$
vanishes on a subset of $\Omega\times\Rm_+$ which is too large in some sense,
then $u(x,t)$ must vanish for all $(x,t)$ in $\Omega\times\Rm_+$.

\section*{Acknowledgments}
The authors are very grateful to Sylvain
Crovisier and Lucien Guillou for fruitful discussions.\\


\end{document}